\documentclass[12pt]{elsart}
\usepackage{amsmath}
\usepackage{amssymb}
\usepackage{amscd}
\usepackage[dvips]{graphics}
\usepackage{amsfonts}

\newcommand{\refm [1]} {(\ref{#1})}

\newcounter{theorem}

\newtheorem{theorem}{Theorem}
\newtheorem{lemma}{Lemma}

\newcommand{\E}{\mathbf{E}}
\newcommand{\la}{\label}

\newcommand{\PP}{\mathbb{P}}

\newcommand{\nicef}{{\bf {\cal F}}}
\newcommand{\en}{\end{equation}}

\def\BP{{\Bbb P}}
\def\PP{{\Bbb P}}
\def\BP{{\Bbb P}}
\def\BR{{\Bbb R}}

\def\E{{\Bbb E}}

\def\BC{{\Bbb C}}

\def\proof{\noindent{\bf Proof\ \ }}
\def\qed{\mbox{\rule{0.5em}{0.5em}}}

\def\la{\label}
\def\plnc{{c_n^{(1)}}}
\def\plnC{{C^{(1)}}}

\def\tilc{{c^{(2)}_n}}
\def\tilC{{C^{(2)}}}
\def\astc{{c^{(3)}_n}}
\def\astC{{C^{(3)}}}

\def\deli{\delta_n^{(i)}}
\def\del1{{\delta_n^{(1)}}}
\def\delh{{\hat{\delta}_n}}

\newcommand{\ignore}[1]{}

\def\BP{{\Bbb P}}
\def\BP{{\Bbb P}}
\def\BR{{\Bbb R}}

\def\BC{{\Bbb C}}
\def\calh{{\cal H}}

\newcommand{\be}{\begin{equation}}
\newcommand{\eu}{\end{equation}}
\newcommand{\ber}{\begin{eqnarray}}
\newcommand{\ena}{\end{eqnarray}}
\newcommand{\nin}{\noindent}
\newcommand{\non}{\nonumber}

\newcommand{\halfe}{\frac{1}{2}}

\begin{document}

\begin{frontmatter}

\title{Meinardus' theorem on weighted partitions: extensions and a probabilistic
proof}
\author{Boris L.  Granovsky}
\address
{ Department of Mathematics, Technion-Israel
Institute of Technology, Haifa, 32000, Israel}
\ead{mar18aa@techunix.technion.ac.il}

\author{Dudley Stark}
\address{School of Mathematical Sciences, Queen Mary,
University of London, London E1 4NS, United Kingdom}
\ead{D.S.Stark@maths.qmul.ac.uk}

\author{Michael Erlihson}
\address{Department of Mathematics, Technion-Israel Institute of
Technology, Haifa, 32000, Israel}
\ead{maerlich@techunix.technion.ac.il}

\begin{keyword}
Combinatorial structures \sep Local limit theorem
\MSC  Primary-60J27 \sep Secondary-60K35
\sep 82C22 \sep 82C26
\end{keyword}
\begin{abstract}
The number $c_n$ of weighted partitions of an integer  $n,$ with
parameters
 (weights) $b_k$, $k\geq 1$, is given by the generating function
relationship $\sum_{n=0}^{\infty}c_nz^n=
\prod_{k=1}^\infty(1-z^k)^{-b_k}$. Meinardus(1954) established his
famous  asymptotic formula for $c_n,$ as $n\to \infty,$ under
three conditions on power and Dirichlet generating functions for
the sequence  $b_k.$ We give a probabilistic proof of Meinardus'
theorem with weakened third condition and extend  the resulting
version of the theorem
 from weighted partitions to other two classic types of decomposable
combinatorial structures, which are called assemblies and
selections.
\end{abstract}

\end{frontmatter}

 \newpage

\section{Summary}
In this paper, we combine Meinardus' approach for deriving the
asymptotic formula for the number  of weighted partitions with the
probabilistic method of Khintchine to develop a unified method of
asymptotic enumeration of three basic types of decomposable
combinatorial structures: multisets, selections and assemblies. As
a byproduct of our approach we weaken one of the three Meinardus
conditions. In accordance with these two objectives, the structure
of the paper is as follows. Section 2 presents Meinardus'
asymptotic formula, the presentation being accompanied by remarks
clarifying the context of the three conditions of Meinardus'
theorem. In section 3 we state our main result   which are
asymptotic formulae for numbers of multisets, selections and
assemblies. Sections 4, 5 and 6 are devoted to the proof of the
main theorem, including  the unified representation of basic
decomposable random structures, which is the core of the
probabilistic method considered. In Section~7 we discuss the
striking similarity between the derived asymptotic formulae.

 \section{Meinardus' theorem}
The Euler type generating function $f^{(1)}$ for the numbers
$\plnc, \ n\ge 1$ of weighted partitions of an integer $n$, with
parameters $b_k\ge 0$, $k\geq 1$ is \be\label{fdef}f^{(1)}(z):=
\sum_{n=0}^{\infty}\plnc z^n=
\prod_{k=1}^\infty(1-z^k)^{-b_k},\quad \vert z \vert<1. \eu In
this setting, $b_k$ is interpreted as a number of types of
summands of size $k$. (For example, one can imagine that coins of
a value $k$ are  distinguished by $b_k$ years of their
production). It is also assumed that in a partition, each summand
of size $k$ belongs to one of the $b_k$ types. In the case $b_k=1$
for all $k\geq 1$, $\plnc$ is the number of standard
(non-weighted) partitions of $n$ (with $c_0=1$), while the case
$b_k=k,\ k\ge 1$ conforms to planar partitions, studied
 by Wright, see \cite{A} and the recent paper \cite{Mu} by
 Mutafchiev.
The study of the asymptotics of the general generating function
(\ref{fdef}) was apparently initiated by Brigham who obtained in
\cite{Br} the asymptotic formula, as $n\to \infty$ for the
logarithm of the function, using the Hardy-Ramanujan asymptotic
technique. Meinardus'  approach (\cite{M}) to the asymptotics of
$\plnc$ is based on considering  two generating series for the
sequence \ $b_k\ge 0,\ k\ge 1$: the Dirichlet series $D(s)$ and
the power series $G(z)$, defined by \be D(s)=\sum_{k=1}^\infty b_k
k^{-s}, \quad s=\sigma+it \label{D},\eu \be G(z)=\sum_{k=1}^\infty
b_kz^k, \quad \vert z\vert<1.\label{Laplace}\eu We note that the
function $f^{(1)}(z)$ converges at the point $\vert z\vert<1$ if
and only if the same is true for the function $G(z)$ (see e.g.
Lemma 1.15 in \cite{Bur}).

Meinardus (\cite{M})  established the following seminal asymptotic
formula for $\plnc$, which is presented in \cite{A}. We denote by
$\Re(\bullet)$ and $\Im(\bullet)$ the real and imaginary parts of
a number.

\begin{theorem}[Meinardus]\label{Meinardus}
Suppose that the parameters $b_k\ge 0, \ k\ge 1$ of  weighted
partitions meet the following three conditions:

\begin{itemize}

\item[(i)] The Dirichlet series (\ref{D}) converges in the
half-plane $\sigma>r>0$ and there is a constant $0<C_0\le 1,$ such
that the function $D(s)$, $s=\sigma+it$, has an analytic
continuation to the half-plane \be\label{Hdef}
\calh=\{s:\sigma\geq-C_0\} \eu on which it is analytic except for
a simple pole at $s=r$ with residue $A$.

\item[(ii)] There is a constant $C_1>0$ such that
\be\label{imagbound} D(s)=O\left(|t|^{C_1}\right),\quad t\to\infty
\eu uniformly in $\sigma\ge -C_0$.

\item[(iii)] There are constants $C_2>0$, $\epsilon>0$ such that
the function \be\label{gdef} g(\tau):=G(\exp(-\tau)),\quad
\tau=\delta+2\pi i\alpha, \ \delta>0,\ \alpha\in R,  \eu satisfies
\be\label{gcond} \Re(g(\tau))-g(\delta)\leq
-C_2\delta^{-\epsilon},\quad \vert\arg \tau\vert>\frac{\pi}{4},
\quad 0\neq|\alpha|\leq 1/2, \eu for $\delta>0$ small enough.
\end{itemize}
Then, as $n\to \infty,$ \be\label{casympmul} \plnc\sim \plnC
n^{\kappa_1}\exp\left(n^{r/(r+1)}\big(1+{1\over r}\big)
\big(A\Gamma(r+1)\zeta(r+1)\big)^{1/(r+1)}\right), \eu where
$$
\kappa_1=\frac{2D(0)-2-r}{2(1+r)}
$$
and
$$
\plnC=e^{D^\prime(0)}
\left(2\pi(1+r)\right)^{-1/2}\left(A\Gamma(r+1)\zeta(r+1)\right)
^{\kappa_2},
$$
where
$$\kappa_2=\frac{1-2D(0)}{2(1+r)}.$$
\end{theorem}

Meinardus also gave a bound on the rate of convergence which
we have omitted in the statement of Theorem~\ref{Meinardus}.

At this point we wish  to make a few clarifying comments on the
three Meinardus conditions $(i), (ii), (iii)$.

\begin{itemize}
\item The Ikehara-Wiener Tauberian theorem on Dirichlet series
cited below tells us that condition $(i)$ implies a bound on
the rate of growth, as $k\to\infty$ of the coefficients  $b_k$ of
the Dirichlet series $D(s)$ in \refm[D].
\begin{theorem}[Wiener-Ikehara](see Theorem 2.2, p.122 in \cite{K})  \label{Wiener}

Suppose that the Dirichlet series $D(s)=\sum_{k=1}^\infty a_k
k^{-s}$ is such that the function $D(s)-\frac{A}{s-1}$ has an
analytic continuation to the closed half-plane $\Re(s)\ge 1$.
Then, \be \sum_{k=1}^na_k\sim An, \quad n\to \infty.\label{win}
\eu
\end{theorem}
We will use the fact that \refm[win] implies \be a_k=o(k),\quad
k\to \infty. \label{ak}\eu To prove this,
 we rewrite \refm[win] as
$${1\over n}\sum_{k=1}^na_k={1\over n} a_n+{1\over n}\sum_{k=1}^{n-1}a_k=A+\epsilon_n,
\quad
\epsilon_n\to 0, \quad n\to \infty,$$ which gives
$${1\over n} a_n+\frac{n-1}{n}(A+\epsilon_{n-1})=A+\epsilon_n.$$
Consequently, $\lim_{n\to\infty}(1/n) a_n=0.$

Now set $a_k=k^{-r+1}b_k, \ k\ge 1,$ where $b_k,\ k\ge 1$ satisfy
Meinardus' conditions (i) and (ii). Since $C_0,r>0$, the sequence
$a_k$ obeys the conditions of the Wiener-Ikehara theorem, so that
we get from \refm[ak] the  bound:\be\label{1bnbound} b_k=o(k^r),
\quad k\to \infty. \eu \item Functions satisfying Meinardus'
condition $(ii)$ are called {\em of finite order} in the
corresponding domain. It is known (see e.g. \cite{Ti}, p. 298)
that the sum $ D$ of a Dirichlet series is a function of a finite
order in the half-plane of the convergence of the series. Thus,
condition $(ii)$ requires that the same holds also for the
analytic continuation of $D$ in the domain $\calh$.

 \item  We show below that the condition $(iii)$ is associated
with bounding the so-called zeta sum known from the theory of the
Riemann Zeta function. In fact,
$$\Re(g(\tau))-g(\delta)=-2\sum_{k=1}^\infty
 {b_k{e^{-k\delta}\sin^2(\pi k\alpha)}}, \quad
 \delta>0,\quad\alpha\in \BR,
 $$
which allows us to reformulate  (\ref{gcond}) as \be\label{gcond2}
2\sum_{k=1}^{\infty }
 {b_k{e^{-k\delta}\sin^2(\pi k\alpha)}}
\geq C_2 \delta^{-\epsilon},\quad
0<{\delta\over{2\pi}}<|\alpha|\leq 1/2, \eu for $\delta>0$ small
enough and some $\epsilon>0.$

The verification of condition $(iii)$ in the forthcoming
Lemma~\ref{polygrowth} relies on the lower bound \refm[sin]
below, for the sum $\sum_{k=1}^P
 \sin^2(\pi k\alpha),\ \alpha\in \BR.$
This bound can be derived from the following bound on the zeta sum
in the left hand side of \refm[bnd]  (see \cite{Kar}, p. 112,
Lemma 1): \be \left| \sum_{k=1}^Pe^{2\pi i k\alpha}\right|\le
\min\left\{P,\,\frac{1}{2\parallel \alpha\parallel}\right\}, \quad
P>1, \quad\alpha\in \BR, \label{bnd}\eu where $\parallel
\alpha\parallel$ denotes the distance from $\alpha$ to the nearest
integer. It follows from (\ref{bnd}) that for all $\alpha\in \BR$
$$
2\sum_{k=1}^P
 \sin^2(\pi k\alpha)
\ge P-\left| \sum_{k=1}^Pe^{2\pi i k\alpha}\right| \ge
P-\min\left\{P,\,\frac{1}{2\parallel \alpha\parallel}\right\},
$$
which is convenient to rewrite as \be \label{sin}
2\sum_{k=1}^P\sin^2(\pi k\alpha)\ge
P\left(1-\min\left\{1,\,\frac{1}{2P\parallel
\alpha\parallel}\right\}\right). \eu Under the assumptions in
Meinardus' condition $(iii)$,
  \be \label{cond} 0\neq\parallel \alpha\parallel=\vert \alpha\vert\le 1/2, \quad
\vert \alpha\vert  \delta^{-1}>\frac{1}{2\pi}.\eu Setting  in
(\ref{sin}) \be P=P(\alpha,\delta)=\left[\frac{1+\vert \alpha
\vert \delta^{-1}}{2\vert \alpha\vert}\right]\ge 1,\la{P}\eu where
$[x]$ denotes the integer part of  $x$ and $\delta>0$ is small
enough, we get the desired bound,

\be\label{bnd0} 2\sum_{k=1}^P\sin^2(\pi k\alpha)\ge
\frac{\delta^{-1}}{2}, \eu provided (\ref{cond}) holds. It follows
from the above that for any fixed $k_0\ge 1$, and any
$0<\epsilon=\epsilon(\delta;k_0)<1/2,$ \be
2\sum_{k=k_0}^P\sin^2(\pi k\alpha)\ge
\left(\frac{1}{2}-\epsilon\right) \delta^{-1}:=c\delta^{-1},
 \label{bnd1}\eu if
$\delta>0$ is small enough and (\ref{cond}) holds.

In the proof of Lemma ~\ref{polygrowth} below we will also use the
fact that under the condition \refm[cond], the choice \refm[P] of
$P$  provides,  \be\label{Pyupper} P\delta\le{1\over
2}\Big(1+\frac{1}{\vert \alpha \vert \delta^{-1}}\Big)<{1\over
2}(1+2\pi):=d. \eu
\end{itemize}

 It seems not to have been noticed
that  Meinardus' condition $(iii)$ is rather easily satisfied,
as is shown in the following lemma.
\begin{lemma}\label{polygrowth}
Let the sequence $\{b_k\}$ be such that  $b_k\geq \rho k^{r-1}, \
k\ge k_0 $ for some $k_0\ge 1$ and  some constants $\rho,r>0.$
Then (\ref{gcond2}) is satisfied.
\end{lemma}
\proof  Because of (\ref{P}), $P\geq \frac{1}{2}\delta^{-1}$ and
therefore $P>k_0$ for $\delta>0$ small enough. We have,
\begin{eqnarray*}
2\sum_{k=1}^{\infty }b_k e^{-k\delta}\sin^2(\pi k\alpha)&\geq&
2\sum_{k=k_0}^P \rho k^{r-1} e^{-k\delta}\sin^2(\pi k\alpha)\\&\geq&
2\rho e^{-P\delta}\sum_{k=k_0}^P k^{r-1}\sin^2(\pi k\alpha):=Q.
\end{eqnarray*}
In order to get the needed lower bound on $Q$ implied by
(\ref{gcond2}), we need to distinguish between the following two
cases: $(i)\ 0<r<1$ and $(ii)\ r\ge 1$. Applying (\ref{bnd1}) and
(\ref{Pyupper}) we have  in  case $(i),$ $$Q\ge \rho
e^{-P\delta}P^{r-1}c\delta^{-1}\ge \rho
e^{-d}(P\delta)^{r-1}c\delta^{-r}\ge \rho
e^{-d}d^{r-1}c\delta^{-r}$$ and in  case $(ii),$
$$Q\ge\rho e^{-d}k_0^{r-1}c\delta^{-1}.$$ Therefore,
(\ref{gcond2}) is satisfied with $\epsilon=r$ in  case $(i)$ and
with $\epsilon=1$ in  case $(ii)$. \hfill\qed \vskip .5cm We note
that in \cite{Mu} the validity of condition $(iii)$ was verified
in the particular case of planar partitions ($b_k=k$,\, $k\ge 1$),
via a complicated analysis  of  the power series expansion of the
function $\Re(g(\tau))-g(\delta).$ \vskip .5cm

\nin {\bf Example 1} Let $b_k=\rho k^{r-1},\ \rho,r>0, \ k\ge 1$.
Such weighted partitions are associated with the generalized
Bose-Einstein model of ideal gas (see \cite{V1}). In this case,
$D(s)=\rho\zeta(s-r+1)$, where $\zeta$ is the Riemann zeta
function. Thus, $D(s)$ has only one simple pole at $s=r>0$ with
the residue $A=\rho$ and it has a meromorphic analytic
continuation to the whole complex plane $\BC$. These facts
together with Lemma~\ref{polygrowth} show that all three
of Meinardus' conditions $(i),(ii),(ii)$ hold. In the case
considered the values $D(0)=\rho\zeta(1-r)$  and
$D^\prime(0)=\rho\zeta^\prime(1-r)$ in the asymptotic formula
\refm[casympmul] can be found explicitly from the functional
relation for the function $\zeta, $ as  is explained in \cite{Mu}.
 In particular, for  standard partitions
($\rho=r=1$),
$$D(0)=\zeta(0)=-\frac{1}{2}, \quad D^\prime(0)=\zeta^\prime(0)=-\frac{1}{2} \log 2\pi,$$
while for  planar partitions ($\rho=1, r=2$),
$$D(0)=\zeta(-1)=-\frac{1}{12}, \quad D^\prime(0)=\zeta^\prime(-1)=
2\int_0^\infty \frac{w \log w}{e^{2\pi w}-1}dw.$$   For an
arbitrary $r>0$, the expressions for $D(0), D^\prime(0)$ include
the integral
$$\int_0^\infty \frac{w^{r-1} \log w}{e^{2\pi w}-1}dw.$$

%
{\bf Example 2} The purpose of this example is to show that
conditions $(i)$ and $(ii)$ of Theorem~\ref{Meinardus} do not
imply condition $(iii)$ in the same theorem. Let
$$ b_k= \left\{
\begin{array}{c c}
1, & {\rm if \ } 4 | k\\
0, & {\rm if \ } 4\!\!\not| k.
\end{array}
\right.
$$
Let $\alpha=1/4$ in the sum $\sum_{k=1}^\infty b_k e^{-\delta
k}\sin^2(\pi k\alpha)$. Then, because for all $k$ either $b_k=0$
or $\sin^2(\pi k/4)=0$,
$$
\sum_{k=1}^\infty b_k e^{-\delta k}\sin^2(\pi k/4)=0 \quad
$$
and therefore (\ref{gcond2}) is not satisfied. However,
$$
D(s)=\sum_{j=1}^\infty (4j)^{-s}=4^{-s}\zeta(s),
$$
which clearly satisfies the first two of Meinardus' conditions
because the function $4^{-s}$ is entire and $\vert
4^{-s}|=4^{-\sigma}\leq 4^{C_0}$ for $s\in\calh$, where $\calh$ is
given by (\ref{Hdef}).
 \section{Statement of the main result}
Our main result, Theorem~\ref{main} below,  achieves two
objectives: weakening the Meinardus condition $(iii)$ and
extending the resulting version of the Meinardus theorem from
weighted partitions to other two  types of classic decomposable
combinatorial structures.

We first recall that a decomposable structure is defined as a
union of indecomposable components of various sizes. It is known
(see \cite{ABT,AT}) that the three types of decomposable
combinatorial structures: multisets, which are also called
weighted partitions, selections and assemblies, encompass the
variety of classic combinatorial objects.  Weighted partitions are
defined as in the previous section, selections are defined as
weighted partitions in  which no component type appears more than
once and assemblies are combinatorial objects composed of
indecomposable components which are formed from labelled elements.
Each decomposable structure is essentially determined by the
number  of types of its indecomposable components having a given
size $k$. We denote this number by $b_k$ for weighted partitions
and selections and by $m_k$ for assemblies. In the case of
assemblies we denote $b_k=m_k/k!$, so that in all three cases
$b_k, \ k\ge 1$ are parameters defining a structure.  In what
follows we will use the notation $\bullet^{(i)}, \ i=1,2,3$ for
quantities related to weighted partitions, selections and
assemblies respectively. Given a sequence $b_k, \ k\ge 1$, we
define $c_n^{(i)}=s_n^{(i)}$ for $i=1,2$ and define
$c_n^{(3)}=s_n^{(3)}/n!$, where $s_n^{(i)}$ denotes in all three
cases  the number of combinatorial structures of type $i$ having
size $n$.

\begin{theorem}\label{main}
Suppose that the parameters $b_k, \ k\ge 1$ meet  Meinardus'
conditions (i) and (ii)  as well as the condition

\begin{itemize}
\item[(iii')] For $\delta>0$ small enough and any
$\epsilon>0,$
$$
2\sum_{k=1}^\infty
 {b_k{e^{-k\delta}\sin^2(\pi k\alpha)}}\ge\left(1+{r\over 2}+\epsilon\right)
M^{(i)} |\log \delta|,$$
$$ \quad  \sqrt{\delta} \leq
|\alpha|\leq 1/2,\quad i=1,2,3,
$$
\end{itemize}
where  the constants $M^{(i)}$ are defined by
$$M^{(i)}=
  \begin{cases}
    {4\over{\log 5}} , & \text{if } \quad i=1,\\
    4 , & \text{if}\quad i=2,\\
     1, & \text{if} \quad i=3.
  \end{cases}
$$

Then the asymptotics for $c_n^{(i)}, \ i=1,2,3,$ \ as \ $n\to
\infty,$ are given respectively by Meinardus' formula
\refm[casympmul], and by the formulae \refm[sl],\refm[as]
below:\be \label{sl}\tilc\sim \tilC n^{-\frac{r+2}{2r+2}}
\exp\left(n^{\frac{r}{r+1}}\big(1+{1\over r}\big)
\big(A(1-2^{-r})\zeta(r+1)\Gamma(r+1)\big)^{\frac{1}{r+1}}\right),
\eu where
$$ \tilC=2^{D(0)} \Big(2\pi(1+r)\Big)^{-1/2}
\Big(A\Gamma(r+1)(1-2^{-r})\zeta(r+1)\Big)^{\frac{1}{2r+2}},$$ and
\be \label{as} \astc\sim \astC n^{-\frac{r+2}{2r+2}}
\exp\left(n^{\frac{r}{r+1}}\big(1+{1\over r}\big)
\big(A\Gamma(r+1)\big)^{\frac{1}{r+1}}\right), \eu where
$$ \astC=e^{D(0)}\Big(2\pi(1+r)\Big)^{-1/2}
\left(A\Gamma(r+1)\right)^{\frac{1}{2r+2}}. $$
\end{theorem}

{ \bf Remark} H-K Hwang (\cite{H}) applied the approach of
Meinardus to the study of the  asymptotics of the number of
summands, say $\omega_n,$ in weighted partitions and
selections, which he called unrestricted and restricted partitions
respectively. In the first case a local limit theorem for a
properly scaled $\omega_n$  was obtained in \cite{H} under the
three conditions of Meinardus. Regarding restricted partitions,
the author claimed the same under Meinardus' conditions $(i)$,
$(ii)$ and a condition similar to our $(iii^\prime),$ but the
proof contains an error in bounding the function $G_\theta(r)$ on
p.109.

 \noindent {\bf Example 3} This   example
 satisfies all three conditions of our Theorem~\ref{main}, but does
not satisfy condition $(iii)$ of Theorem~\ref{Meinardus}. Let
$b_k$, $k\geq 1$, be defined by
$$
b_k= \left\{
\begin{array}{l l}
12e^7\big(\frac{\log k}{k}\big) & {\rm if \ } 4\!\!\not| k,\\
12e^7\big(50+\frac{\log k}{k}-2\frac{\log(k/4)}{k/4}\big)
& {\rm if \ } 4 | k {\rm \ and \ } 16\!\!\not| k,\\
12e^7\big( 50+\frac{\log
k}{k}-2\frac{\log(k/4)}{k/4}+\frac{\log(k/16)}{k/16} \big) & {\rm
if \ } 16 | k.
\end{array}
\right.
$$
Note that because $0\leq \frac{\log x}{x}\leq e^{-1}$ for $x\geq
1,$ it follows that $b_k\geq 12e^7\big(\frac{\log k}{k}\big),$ for
all $k\ge 1.$ The Dirichlet series $D(s),\ s=\sigma+it,$   for
this choice of $b_k$ converges absolutely for $\sigma>1$ and in
this domain \be\label{Dir}
D(s)=12e^7\Big(-(1-4^{-s})^2\zeta^\prime(s+1)+50\cdot
4^{-s}\zeta(s)\Big), \eu where we have used the fact that
$\zeta^\prime(s+1)=-\sum_{k=1}^\infty \frac{\log k}{k} k^{-s}$. It
is well known that the function $\zeta(s)$ has a simple pole at
$s=1$ and that the Laurent expansion of $\zeta(s+1)$ around $s=0$
is \be\label{zetaexp} \zeta(s+1)={1\over s}+\gamma+\cdots, \eu
where
$\gamma$ is Euler's constant. It follows from (\ref{zetaexp}) that
the function
 $\zeta^\prime(s+1)$ has a unique pole at $s=0$ of order $2$.
As a result, we derive that in \refm[Dir] the first term in the
parentheses is analytic in the whole complex plane $\BC$, while
the function $D$ in (\ref{Dir}) is analytic in $\BC$ except a
simple pole at $s=1$. It is also a known fact that the functions
$\zeta, \zeta^\prime$ satisfy (\ref{imagbound}), from which we
conclude that the same is true for the function $D$ given by
(\ref{Dir}). To show that condition $(iii^\prime)$ of
Theorem~\ref{main} is satisfied, we note that, if $\delta>0$ and
$\epsilon>0$ are small enough then

\begin{eqnarray}
\sum_{k=1}^\infty b_k e^{-k\delta}\sin^2(\pi k\alpha)&\geq&
\sum_{k=1}^P
b_k e^{-k\delta}\sin^2(\pi k\alpha)\nonumber\\
&\geq& 12e^7 e^{-d} \sum_{k=1}^P \frac{\log k}{k}
\sin^2(\pi k\alpha)\nonumber\\
&\geq& 12e^7 e^{-d}\frac{\log P}{P} \sum_{k=3}^P
\sin^2(\pi k\alpha)\nonumber\\
&\geq& 12e^{7-d}\frac{\log (d \delta^{-1})}{d \delta^{-1}}
(\frac{1}{4}-\epsilon)\delta^{-1}\label{examstep}\\
&>& 6\log( \delta^{-1})\label{lowersum}, \ \
(2\pi)^{-1}\delta<|\alpha|\leq 1/2,
\end{eqnarray} where
we have used (\ref{bnd1}) and (\ref{Pyupper}) at (\ref{examstep})
and the fact that $3.5<d<4$ in the last step. Since in the case
considered $r=1$, the condition $(iii^\prime)$ is indeed satisfied
for all three types of random structures. Finally, to show that
condition $(iii)$ of Theorem~\ref{Meinardus} in the form
(\ref{gcond2}) is not satisfied, we set $\alpha=1/4$ in the left
hand side of (\ref{gcond2}) to obtain for $ \delta\to 0,$
\begin{eqnarray*}\sum_{k=1}^\infty b_k e^{-k\delta}\sin^2(\pi
\frac{k}{4})\le 12e^7\sum_{k=1}^\infty \frac{\log k}{k}
e^{-k\delta}
\sim  12e^7\int_1^\infty \frac{\log x}{x} e^{-x\delta}dx  = \nonumber\\
12e^7\int_{\delta}^1 \frac{\log (\delta^{-1}x)}{x} e^{-x}dx+
12e^7\int_1^\infty \frac{\log (\delta^{-1}x)}{x} e^{-x}dx=
O\left(\log^2 (\delta^{-1})\right).\end{eqnarray*}

{\bf Example 4} Consider the assembly of forests, for which
components consist of labelled linear trees. The number of such
components on $k$ vertices is $m_k=k!$ and so $b_k=1$, just as for
integer partitions. The asymptotic number of labelled linear
forests is thereby given by \refm[as] in Theorem~\ref{main} with
$r=1$, $A=1$. We note that the number of labelled linear forests
on $n$ vertices equals the number of path coverings of a
complete graph on $n$ vertices.

\section{A unified probabilistic representation for decomposable combinatorial structures.}

It has been recently understood (see \cite{Pi,V1}) that the three
main types of decomposable random structures: assemblies,
multisets and selections, are induced by a class of probability
measures on the set of integer partitions, having a multiplicative
form. Vershik (\cite{V1}) calls the measures multiplicative, while
Pitman (\cite{Pi},\cite{BP}) refers to them as Gibbs partitions.
Equivalently, in combinatorics it is common to view the structures
above as the ones generated by the conditioning relation (see
{\cite{ABT}) or by the Kolchin generalized allocation scheme
(\cite{kol}). Our asymptotic analysis is based on the unified
Khintchine type probabilistic representation of the number of
decomposable structures of size $n$. Recall that  we agree
that the number of non labelled structures (weighted partitions
and selections) is denoted  by $c_n^{(1)}$ and $c_n^{(2)}$
respectively, and  the number of labelled structures (assemblies)
by $n!c_n^{(3)}$. In all three cases the probabilistic
representation of $c_n$ is constructed as follows. Let $f$ be a
generating function of a sequence $\{c_n\}$ associated with some
decomposable structure:
$$
f(z)=\sum_{n\ge 1}c_n z^n.
$$
A specific feature of decomposable structures is that the
generating function $f$ has the following multiplicative form: $$
f=\prod_{k\ge 1}S_k,$$  where $S_k$ is a generating function for
some nonnegative sequence $\{d_k(j),\ j\ge 0,\ k\ge1\},$ i.e.\be
S_k (z)=\sum_{j\ge 0}d_k(j)z^{kj}, \ k\ge 1.\la{g3}\eu We now set
$z=e^{-\delta +2\pi i\alpha},\ \alpha\in [0,1]$ and use the
orthogonality property of the functions $e^{-2\pi i\alpha n}, \
n\ge 1$ on the set $\alpha\in[0,1]$, to get
\begin{eqnarray}
c_n& =& e^{n\delta} \int_0^1 f\left(e^{-\delta+2\pi
i\alpha}\right)
   e^{-2\pi i\alpha n} d\alpha \nonumber\\ &=& e^{n\delta} \int_0^1
 \prod_{k=1}^n\Big( S_k\left(e^{-\delta+2\pi i\alpha}\right)\Big) e^{-2\pi i\alpha
n} d\alpha, \quad n\ge 1, \nonumber \la{qe}
\end{eqnarray}
where $\delta$ is a free parameter.\ We denote by \be
f_n:=\prod_{k=1}^n S_k, \ \ n\ge 1, \la{g2} \eu  the truncated
generating function. Next, we attribute a probabilistic meaning to
the expression in the right hand side of (\ref{qe}) by defining
the independent integer valued random variables $Y_k, \ k\ge 1:$
\be \PP(Y_k=jk)= \frac{d_k(j) e^{-\delta kj}}{S_k(e^{-\delta })},
\quad j\ge 0,\quad k\ge 1.\la{Y}\eu and observing that \be
\phi_n(\alpha):=\prod_{k=1}^n \frac{S_k(e^{-\delta+2\pi
i\alpha})}{S_k(e^{-\delta})}=\frac{f_n(e^{-\delta+2\pi
i\alpha})}{f_n(e^{-\delta})}=\E \Big(e ^{2\pi i\alpha Z_n}\Big),
\quad \alpha\in R \la{ch} \eu is the characteristic function of
the random variable \be Z_n:=\sum_{k=1}^{n}{Y_k}. \la{zdef}\eu We
have arrived at the desired representation:
 \be
c_n=e^{n\delta} f_n(e^{-\delta})\PP\left(Z_n=n\right),\quad n\ge
1.
 \la{rep}\eu

In accordance with the principle of the probabilistic method
considered, we will
 choose in (\ref{rep}) the free parameter $\delta=\delta_n$ to be the
 solution of the equation
 \ber
 \E{Z_n}=n, \quad n\ge 1,
 \la{deltandef}
 \ena
after we show in the next section that for the three classic
combinatorial structures the solution to (\ref{deltandef}) exists
and is unique.

It can be easily seen from (\ref{g3}), (\ref{g2}), (\ref{Y})
 that \be \E{Z_n}=(\E{Z_n})(\delta)=-\left(\log
f_n(e^{-\delta})\right)^\prime, \quad \delta>0. \la{EZ}\eu It is
interesting to note that in the context of thermodynamics, the
quantity $\log f_n(e^{-\delta})$ has a meaning of the  entropy of
a system. This important fact that clarifies the choice of the
free parameter  was observed already by Khintchine
(\cite{Kh},Chapter VI), in the course of his study of classic
models of thermodynamics.

 From this point on, our study will be restricted to the three
above mentioned classic combinatorial structures: multisets
(weighted partitions), selections and assemblies. Recalling the
forms of their generating functions $f^{(i)}, \ i=1,2,3$ (see
\cite{ABT}) and denoting
$\nicef^{(i)}(\delta)=f^{(i)}\left(e^{-\delta}\right),\ \delta>0,\
i=1,2,3$, we obtain
\begin{eqnarray}
\nicef^{(1)}(\delta)&=&\prod_{k\ge 1}{(1-e^{-k\delta})^{-b_k}},\nonumber\\
\nicef^{(2)}(\delta)&=&\prod_{k\ge 1}{(1+e^{-k\delta})^{b_k}},\nonumber\\
\nicef^{(3)}(\delta)&=&\exp\left(\sum_{k\ge
1}b_ke^{-k\delta}\right).\la{F}
\end{eqnarray}

Now it is easy to derive from \refm[Y] and \refm[F] that  the
following three types of distributions for the random variables
$\frac{1}{k}Y_k$ in \refm[Y]:
 Negative
 Binomial $\left(b_k;e^{-\delta k}\right),$
Binomial $\left(b_k;\frac{e^{-\delta k}}{1+e^{-\delta k}}\right)$
and Poisson $\left(b_ke^{-\delta k}\right)$, produce
respectively
 $c_n^{(i)}, \ i=1,2,3$ in  the representation  \refm[rep].

The representation (\ref{rep}) for assemblies was  obtained in
 \cite{frgr1}, while the one for general multisets and selections
 was obtained in \cite{GS}.
The corresponding truncated generating functions $f_n^{(i)}(z)$
are
{\allowdisplaybreaks
\begin{eqnarray}
f_n^{(1)}(z)&=&\prod_{k=1}^n{(1-z)^{-b_k}},\nonumber\\
f_n^{(2)}(z)&=&\prod_{k=1}^n{(1+z)^{b_k}},\nonumber\\
f_n^{(3)}(z)&=&\exp\left(\sum_{k=1}^nb_kz^k\right).\la{f}
\end{eqnarray}
}

 Consequently, in the three cases considered
 the
equation (\ref{deltandef}) takes the forms (\ref{e1})-(\ref{e3})
derived from \refm[EZ]
{\allowdisplaybreaks
\begin{eqnarray}
\sum_{k=1}^{n}{\frac{kb_ke^{-k\delta^{(1)}_n}}{1-e^{-k\delta^{(1)}_n}}}
&=&n, \la{e1}\\
 \sum_{k=1}^{n}{\frac{kb_ke^{-\delta^{(2)}_n
k}}{1+e^{-\delta^{(2)}_n k}}}&=&n,  \la{e2}\\
\sum_{k=1}^{n} {kb_ke^{-\delta^{(3)}_n k}}&=&n. \la{e3}
\end{eqnarray}
}

\section{Preliminary asymptotic results}
In this section we find asymptotics for solutions to
(\ref{e1})-(\ref{e3}).
 \begin{lemma}\label{prodasymp}
Suppose that the sequence $b_k\ge 0,\ k\ge 1$ is such that  the
associated Dirichlet generating function $D$  satisfies the
conditions $(i)$ and $(ii)$ of Theorem~\ref{Meinardus}.
 Then

\nin (i) As $\delta\to 0^+,$
\begin{eqnarray}
\nicef^{(1)}(\delta)&=&\exp{\left(A\Gamma(r)\zeta(r+1)\delta^{-r}
  -D(0)\log{\delta}+D'(0)+O\left(\delta^{C_0}\right)\right)},\label{prod1}\\
\nicef^{(2)}(\delta)&=&\exp{\left(A\Gamma(r)(1-2^{-r})\zeta(r+1)
 \delta^{-r}+D(0)\log 2+O\left(\delta^{C_0}\right)\right)},
\label{prod2}\\
\nicef^{(3)}(\delta)&=& \exp{\left(A\Gamma(r)\delta^{-r}+D(0)+
O\left(\delta^{C_0}\right)\right)}, \label{prod3} \end{eqnarray}

\nin whereas  asymptotic expressions for the derivatives
$$\Big(\log\nicef^{(i)}(\delta)\Big)^{(k)},\ \ i=1,2,3;\ k=1,2,3$$ are given by the formal differentiation of the logarithms of
(\ref{prod1})-(\ref{prod3}):
\begin{eqnarray}
\Big(\log\nicef^{(1)}(\delta)\Big)^{(k)} &=&(-1)^k
A\Gamma(r+k)\zeta(r+1)\delta^{-r-k}
  +(-1)^{k} (k-1) D(0)\delta^{-k}+O\left(\delta^{C_0-k}\right), \nonumber\\
\Big(\log\nicef^{(2)}(\delta)\Big)^{(k)} &=&(-1)^k
A\Gamma(r+k)(1-2^{-r})\zeta(r+1)
 \delta^{-r-k}+O\left(\delta^{C_0-k}\right), \nonumber\\
\Big(\log\nicef^{(3)}(\delta))\Big)^{(k)}&=&
(-1)^kA\Gamma(r+k)\delta^{-r-k}+O\left(\delta^{C_0-k}\right).
\label{prime}
\end{eqnarray}

 \nin (ii) Each of the equations (\ref{e1})-(\ref{e3}) has a
unique solution $\delta_n^{(i)}$ such that \\ $\delta^{(i)}_n\to
0, \ n\to \infty, \ i=1,2,3.$

Moreover,

\nin (iii) As $n\to \infty,$

$$
\delta_n^{(1)}=
\left(A\Gamma(r+1)\zeta(r+1)\right)^{\frac{1}{r+1}}
n^{-\frac{1}{r+1}}+\frac{D(0)}{r+1}n^{-1}+O(n^{-1-\beta}),
$$
where
\be
\beta=
  \begin{cases}
\label{beta}
    \frac{C_0}{r+1} , & \text{if}\ \ r\ge C_0 \\
    \frac{r}{r+1} , & \text{otherwise};
  \end{cases}
\eu
\begin{eqnarray} \delta_n^{(2)}=
\left(A\Gamma(r+1)(1-2^{-r})\zeta(r+1)\right)^{{1\over{r+1}}}n^{-\frac{1}{r+1}}&+&
O\left(n^{-1-\beta}\right),\quad \nonumber\\
{\rm where}\ \beta=\frac{C_0}{r+1}; \la{delta2}
\end{eqnarray}
\begin{eqnarray}
\delta_n^{(3)}=
\left(A\Gamma(r+1)\right)^{{1\over{r+1}}}n^{-\frac{1}{r+1}}&+&
O\left(n^{-1-\beta}\right),\nonumber\\ {\rm where \
}\beta=\frac{C_0}{r+1}. \la{delta3}
\end{eqnarray}

\nin (iv) As $n\to \infty$, $f^{(i)}_n(e^{-\deli}),$ $i=1,2,3$
have the asymptotic expansions of the right hand sides of
(\ref{prod1})- (\ref{prime}) respectively, with $\delta=\deli$.
\end{lemma}

 \proof (i) First consider the case of weighted partitions.
Following the Meinardus approach, we will use the fact that
$e^{-u}$, $u>0$, is the Mellin transform of the Gamma function:
\be e^{-u}=\frac{1}{2\pi i}\int_{v-i\infty}^{v+i\infty}
u^{-s}\Gamma(s)\,ds,\quad u>0, v>0. \la{Mellin} \en Expanding
$\log ~\nicef^{(1)}(\delta)$ in (\ref{F}) as
$$
\log
~\nicef^{(1)}(\delta)=-\sum_{k\ge 1}b_k\log(1-e^{-\delta k}) =\sum
_{j\ge 1}\frac{1}{j}\sum_{k\ge 1}b_ke^{-\delta kj}
$$
and substituting (\ref{Mellin}) with $v=1+r$ gives
\be\label{intrep} \log ~\nicef^{(1)}(\delta)=\frac{1}{2\pi
i}\int_{1+r-i\infty}^{1+r+i\infty}
\delta^{-s}\Gamma(s)\zeta(s+1)D(s)ds. \en By Meinardus'
condition (i), the function $D$ has a simple pole at $r>0$ with
residue $A$, which says that the integrand in (\ref{intrep}) has a
simple pole at $s=r$ with  residue
$A\delta^{-r}\Gamma(r)\zeta(r+1)$.
 Next,
from the Laurent expansions at $s=0$ of the Riemann Zeta
 function $\zeta(s+1)=\frac{1}{s}+\gamma+\ldots$ and
the Gamma function $\Gamma(s)=\frac{1}{s}-\gamma+\ldots,$ where
$\gamma$ is
 Euler's constant, and the Taylor series expansions at $s=0$ of the
two remaining factors of the integrand in \refm[intrep], one
concludes that the integrand has also a pole of a second order at
$s=0$ with residue $D^\prime(0)-D(0)\log\delta$. We also recall
that the only poles of $\Gamma(s)$ are at $s=-k, \ k=0,1,\ldots$
Hence, in the complex domain $ -C_0\le \Re( s)\le 1+r,$ with $0\le
C_0<1,$ the integrand has only two poles at $0 $ and $r$ with the
above residuals. We now apply the residue theorem for the
integrand in (\ref{intrep}), over the above domain. The assumption
(\ref{imagbound}) and the following two properties of Zeta and
Gamma functions
$$
\zeta(\sigma+1+it)=O\left(\vert t\vert^{C_2}\right), \ t\to
\infty, \quad C_2>0,
$$
$$
\Gamma(\sigma+it)=O\left(\vert
t\vert^{C_3}\exp(-\frac{\pi}{2}\vert t\vert)\right), \ t\to
\infty, \quad C_3>0,
$$
uniformly in $\sigma$, allow us to conclude that  the  integral
of the integrand considered, over the horizontal contour
$-C_0\le\Re(s)\le 1+r$, $\Im(s)=t$, tends to zero, as $t\to \infty,$ for
any fixed $\delta.$ Thus, we are able to rewrite (\ref{intrep}) as
\begin{eqnarray}\label{expand}
\log ~\nicef^{(1)}(\delta)
&=&A\delta^{-r}\Gamma(r)\zeta(r+1)-D(0)\log \delta + D^\prime(0)\nonumber\\
&&+\,\frac{1}{2\pi
i}\int_{-C_0-i\infty}^{-C_0+i\infty}\delta^{-s}\Gamma(s)\zeta(s+1)D(s)ds.
\end{eqnarray}
Moreover, the previous two bounds and  the bound \refm[imagbound]
in Meinardus' condition (ii)  imply that the integral in
(\ref{expand}) is bounded by
\begin{eqnarray*}
&&
\left|\frac{1}{2\pi i}\int_{-C_0-i\infty}^{-C_0+i\infty}
\delta^{-s}\Gamma(s)\zeta(s+1)D(s)\,ds\right|\\
&=&
O\left(\delta^{C_0}\int_{-\infty}^\infty\exp\left(-\frac{\pi}{2}\vert
t\vert\right) \vert t \vert^{C_1+C_2+C_3}dt\right)\\
&=& O\Big(\delta^{C_0}\Big), \quad \delta\to 0.
\end{eqnarray*}
This proves (\ref{prod1}).

To prove the asymptotic formula for the first derivative
$\Big(\log\nicef^{(1)}(\delta)\Big)^{(1)}$, one has to
differentiate (\ref{expand}) with respect to $\delta$ and then to
estimate the resulting integral in the same way as above.
Subsequent differentiations produce the asymptotic formulae for
$\Big(\log\nicef^{(1)}(\delta)\Big)^{(k)},\ k=2,3$.

The proof of part (i) of the theorem for selections and assemblies
is done in a similar way we now briefly describe. Following
\refm[F], the representation \refm[intrep] conforms to \be
\label{intrep2} \log ~\nicef^{(2)}(\delta)=\frac{1}{2\pi
i}\int_{1+r-i\infty}^{1+r+i\infty}
\delta^{-s}\Gamma(s)(1-2^{-s})\zeta(s+1)D(s)ds\eu and
 \be\label{intrep3}
\log ~\nicef^{(3)}(\delta)=\frac{1}{2\pi
i}\int_{1+r-i\infty}^{1+r+i\infty} \delta^{-s}\Gamma(s)D(s)ds, \eu
for all $\delta>0.$ Accordingly,  the  integrand in \refm[intrep2]
has  a simple pole at $s=r>0$ with  residue
$A\delta^{-r}\Gamma(r)(1-2^{-r})\zeta(r+1)$, and a simple pole  at
$s=0$ with residue $D(0)\log 2$, while the integrand in
\refm[intrep3] has two simple poles at $s=r>0$ and $s=0$ with
residues $A\delta^{-r}\Gamma(r)$ and $D(0)$ respectively. As a
result, we obtain \refm[prod2] and \refm[prod3].

(ii) We see that the left hand sides of the equations
(\ref{e1}-\ref{e3}) are decreasing as $\delta\ge 0$ in such a way
that for a fixed $n,$
 in all the three cases the left hand sides tend to
$0$ as $\delta\to +\infty,$ while as $\delta\to 0$ the left hand
sides tend to $+\infty$, ${1\over 2}\sum_{k=1}^n kb_k$ and
$\sum_{k=1}^n kb_k$ respectively. We now make use of
Theorem~\ref{Wiener} to get a lower bound \refm[kbn] below on the
sum $\sum_{k=1}^n kb_k$ when the sequence $\{b_k\}$ obeys
Meinardus' conditions (i) and (ii). We set $a_k:=k^{-r+1}b_k, \
k\ge 1$ and let $\tilde D(s)$ denote the Dirichlet series
$\tilde{D}(s)=\sum_{k\geq 1} a_k k^{-s}$.  Since $C_0,r>0$, the
function $\tilde D$ satisfies the conditions of Wiener-Ikehara
theorem, with the constant $A$ as in Meinardus' condition (i).
Consequently,
$$
\sum_{k=1}^nk^{-r+1}b_k=\sum_{k=1}^n\frac{kb_k}{k^r}\sim An, \quad
n\to \infty,
$$
from which it follows that for sufficiently large $n,$ \be
\sum_{k=1}^n kb_k\ge Bn,\la{kbn}\end{equation} for some $B>1.$
This can be easily seen from the bound
$$\sum_{k=1}^n\frac{kb_k}{k^r}\le
\sum_{k=1}^{L-1}\frac{kb_k}{k^r}+ \frac{1}{L^r}\sum_{k=L}^nkb_k$$
with $L<n$ such that $L^rA>1$).
 Moreover, \refm[1bnbound]
implies that  the series
$\lim_{n\to\infty}\E{Z_n^{(i)}}, \ i=1,2,3$ converge for any positive $\delta.$\\
Combining the above facts, we conclude that each of the equations
(\ref{e1}-\ref{e3}) has a unique solution for sufficiently large
$n$ and that the solutions $\delta_n^{(i)}\to 0, \ n\to \infty,\
i=1,2,3.$

(iii) We firstly show that in all three cases,
$$\E{Z_n^{(i)}}= \left(-\log
\nicef^{(i)}(e^{-\delta})\right)^\prime\Big|_{\delta=\delta_n^{(i)}}
+\epsilon^{(i)}(n), \quad \epsilon^{(i)}(n)\to 0,$$\be \quad  n\to
\infty, \quad i=1,2,3.\label{EE}\eu

In the case of weighted partitions, setting
$\delh=n^{-\frac{r+2}{2(r+1)}}$ gives for sufficiently large $n$
\begin{eqnarray}
\sum_{k=n+1}^{\infty}{\frac{kb_ke^{-k\delh}}{1-e^{-k\delh}}}
=O\left(\sum_{k=n+1}^\infty kb_ke^{-k\delh}\right)=
O\left(\sum_{k=n+1}^\infty k^{r+1}e^{-k\delh}\right)&=&\nonumber\\
O\Big(\int_{n+1}^\infty x^{r+1}e^{-x\delh}dx\Big)\to 0, \ n\to
\infty,\label{eng}
\end{eqnarray}
where we have employed (\ref{1bnbound}) and the fact that
$n\delh\to \infty, \ n\to \infty.$ From (\ref{prime}) with $k=1$
we deduce that that $\delta_n^{(1)}>\delh$ for large enough $n,$
which implies that \refm[eng] is valid with $\delh$ replaced by
$\delta_n^{(1)}.$ This  proves \refm[EE] for the case considered.
Consequently, the equation  (\ref{e1}) can be rewritten as

\be
A\Gamma(r+1)\zeta(r+1)
(\delta^{(1)}_n)^{-r-1}+D(0)(\delta^{(1)}_n)^{-1}
+O\left((\delta^{(1)}_n)^{C_0-1}\right)+\epsilon^{(1)}(n)=n, $$$$
\quad \epsilon^{(1)}(n)\to 0, \quad n\to \infty. \la{ase}\eu

We outline here the method of solution for asymptotic equations of
the form (\ref{ase}) common in applications of Khintchine's
method.  Denoting the constant coefficient
$h:=A\Gamma(r+1)\zeta(r+1),$ (\ref{ase}) implies \be
h+D(0)(\delta^{(1)}_n)^{r}+
O\left((\delta^{(1)}_n)^{r+C_0}\right)+o\left((\delta^{(1)}_n)^{r+1}\right)=
n\left(\delta^{(1)}_n\right)^{r+1}. \label{sx} \eu Since
$\delta^{(1)}_n\to 0, \ n\to \infty,$ we obtain from \refm[sx]
that $\delta^{(1)}_n\sim h^{\frac{1}{r+1}} n^{-{\frac{1}{r+1}}}, \
n\to \infty.$ Based on this fact and  the fact that $0<C_0 <1,$ we
get
\begin{eqnarray*}
\delta^{(1)}_n&=&h^{\frac{1}{r+1}}
n^{-\frac{1}{r+1}}+\frac{D(0)}{r+1}n^{-1}+
O\left(\left(\delta^{(1)}_n\right)^rn^{-1}\right)+O\left((\delta^{(1)}_n)^{C_0}n^{-1}\right)
\\
&=& h^{\frac{1}{r+1}}
n^{-{\frac{1}{r+1}}}+\frac{D(0)}{r+1}n^{-1}+O(n^{-1-\beta}),
\end{eqnarray*}
where $\beta$ is as in \refm[beta].
 For selections and assemblies
the analogs of \refm[sx] will be respectively $$
h+O\left(\left(\delta^{(2)}_n\right)^{C_0+r}\right)+
o\left(\left(\delta^{(2)}_n\right)^{r+1}\right)
=n\left(\delta^{(2)}_n\right)^{r+1},$$
$$h=A\Gamma(r+1)(1-2^{-r})\zeta(r+1)$$

and $$ h + O\left(\left(\delta^{(3)}_n\right)^{C_0
  +r}\right)+ o\left(\left(\delta^{(3)}_n\right)^{
  r+1}\right)=n\left(\delta^{(3)}_n\right)^{r+1},$$
$$h= A\Gamma(r+1).$$
Now the same reasoning as for weighted partitions leads to the
solutions \refm[delta2], \refm[delta3].

(iv) In the case of weighted partitions, we have
$$
\log f^{(1)}_n(e^{-\delta^{(1)}_n}) = \log
\nicef_1(\delta^{(1)}_n)+ \sum_{k\ge n+1}b_k\log
(1-e^{-k\delta^{(1)}_n}),
$$
where
$$
\left|\sum_{k\ge n+1}b_k\log (1-e^{-k\delta^{(1)}_n})\right|=
O\left(\sum_{k\ge n+1}b_ke^{-k\delta^{(1)}_n}\right)=
o\left(n^{r+1}\exp\left(-n^{\frac{r}{2r+2}}\right)\right),$$$$
\quad n\to \infty,
$$
by the argument giving (\ref{eng}). The proof of the remaining
parts of the assertion (iv) is similar. \hfill\qed

{\bf Remark} As we mentioned before, Meinardus' proof (see
\cite{A}) of Theorem~\ref{Meinardus} relies on application of the
saddle point method. In accordance with the principle of the
method, the value in question $c_n^{(1)}$ is expressed as \be
c_n^{(1)}=\frac{1}{2\pi i}
\int_{-1/2}^{1/2}\nicef^{(1)}(\tau)e^{n\delta+2\pi
in\alpha}d\alpha, \quad \tau=\delta+2\pi i\alpha,\la{mn}\eu by
virtue of the Cauchy integral theorem. Here the free parameter
$\delta$ is chosen as the minimal value of the function
$\exp{\left(A\Gamma(r)\zeta(r+1)\delta^{-r}+n\delta\right)}$
viewed as an approximation of the absolute value of the  integrand
in (\ref{mn}) at $\alpha=0$. This gives
$\delta=h^{\frac{1}{r+1}}n^{-\frac{1}{r+1}},\
h=A\Gamma(r+1)\zeta(r+1)$ which is the principal term of the
solution $\delta_n^{(1)}$ of (\ref{e1}). It can be seen that,
stemming from this choice of the free parameter, the subsequent
steps of Meinardus' proof are considerably more complicated
compared with ours. Also note that our  choice of the free
parameter is in the core of our  ability to weaken
Meinardus' condition $(iii)$. \vskip .5cm Our next assertion
reveals that the function $\Re(g(\tau))-g(\delta)$ in the left
hand side of the Meinardus' condition $(iii)$ is inherent in the
employed probabilistic method: the function provides an upper
bound for the rate of exponential decay of the absolute value of
the characteristic function $\phi_n$ in \refm[ch], as $n\to
\infty,$ for all three types of random structures considered. The
bounds obtained in the forthcoming lemma are used in the the proof
of our local limit theorem, Theorem~\ref{loc}.

Recall that the function $g(\tau)$ is defined by (\ref{Laplace})
and (\ref{gdef}).
\begin{lemma}\label{gtau} Denote
\begin{eqnarray*}
V(\alpha)=V(\alpha;\delta)&=&\Re(g(\tau))-g(\delta), \quad \tau=\delta+2\pi i\alpha\\
&=&-2\sum_{k=1}^\infty b_ke^{-k\delta}\sin^2(\pi\alpha k), \quad
\delta>0,\quad \alpha\in\BR
\end{eqnarray*}
and let $\delta= \deli, \ i=1,2,3$ be the unique solutions of the
equations (\ref{e1})-(\ref{e3}) respectively. Then, for all
$\alpha\in \BR,$
$$\vert\phi_n^{(i)}(\alpha)\vert
\le (1+\epsilon_n^{(i)})\exp\Big({{V^{(i)}(\alpha)}\over
M^{(i)}}\Big), \quad \epsilon_n^{(i)}=\epsilon_n^{(i)}(\alpha)\to
0, \quad n\to \infty, \ i=1,2,3,$$ where
$V^{(i)}(\alpha)=V(\alpha;\delta_n^{(i)})$ and the constants
$M^{(i)}, \ i=1,2,3$ are as in condition $(iii^\prime)$ of
Theorem~\ref{main}.
\end{lemma}
{\bf Proof} From (\ref{ch}) we have for $n\ge 1$ and $\delta>0$
fixed,
$$\vert \phi_n^{(i)}(\alpha)\vert=\exp{\bigg(\Re\left(\log f^{(i)}_n(e^{-\tau})\right)-
\log
f^{(i)}_n(e^{-\delta})}\bigg):=e^{V_n^{(i)}(\alpha;\delta)},$$
\be \tau=\delta+2\pi i \alpha.\label{vn}\eu Using  (\ref{f}) we
now find bounds for $V_n^{(i)}(\alpha;\delta)$ expressed via
$V^{(i)}(\alpha;\delta)$ in the three cases considered. By the
definition of $V_n^{(i)}(\alpha;\delta)$ as given in \refm[vn],
we have for weighted partitions,

{\allowdisplaybreaks
\begin{eqnarray*}
V_n^{(1)}(\alpha;\delta)&=&\Re\left(-\sum_{k=1}^n
b_k\log\left(\frac{1-e^{-\tau k}}{1-e^{-\delta k}}\right)\right)\\
&=& -{1\over 2}\sum_{k=1}^n b_k\log\left(\frac{1-2e^{-k\delta}\cos(2\pi
\alpha
k)+e^{-2k\delta}}{(1-e^{-\delta k})^2}\right)\\
&=& -{1\over 2}\sum_{k=1}^n b_k\log\Big(1+\frac{4e^{-\delta
k}\sin^2(\pi
\alpha k)}{(1-e^{-\delta k})^2}\Big)\\
&\leq&-{1\over 2}\sum_{k=1}^n b_k\log\Big(1+4e^{-\delta
k}\sin^2(\pi
\alpha k)\Big)\\
&\leq&-{{\log 5}\over 2}\sum_{k=1}^n b_ke^{-\delta k}\sin^2(\pi
\alpha k),\quad \delta >0,\quad\alpha\in\BR,
\end{eqnarray*}
}
where the last inequality is due to the fact that $\log(1+x)\ge
(\frac{\log5}{4})x, \ 0\le x\le 4.$

For selections, in a similar manner,
{\allowdisplaybreaks
\begin{eqnarray*}
V_n^{(2)}(\alpha;\delta)&=&\Re\left(\sum_{k=1}^n
b_k\log\left(\frac{1+e^{-\tau k}}{1+e^{-\delta
k}}\right)\right)\\
&=& {1\over 2}\sum_{k=1}^n b_k\log\left(\frac{1+2e^{-k\delta}\cos(2\pi
\alpha
k)+e^{-2k\delta}}{(1+e^{-\delta k})^2}\right)\\
&=&{1\over 2}\sum_{k=1}^n
b_k\log\Big(1-\frac{4e^{-k\delta}\sin^2(\pi \alpha k)}{(1+e^{-\delta
k})^2}\Big)\\
&\leq&-{1\over 2}\sum_{k=1}^n b_k\frac{4e^{-k\delta}\sin^2(\pi
\alpha k)}{(1+e^{-\delta k})^2}\\
&\leq& -{1\over 2}\sum_{k=1}^n b_ke^{-k\delta}\sin^2(\pi \alpha k),
\quad \delta>0, \quad \alpha\in\BR
\end{eqnarray*}
}

Here the first
inequality is due to the fact that $-\log(1-x)\geq x, \quad 0\le
x\le 1.$ For assemblies, we get straightforwardly
$$V_n^{(3)}(\alpha;\delta)=-2\sum_{k=1}^n
b_ke^{-k\delta}\sin^2(\pi \alpha k), \quad \delta>0, \quad \alpha\in
R.$$

Finally, setting $\delta=\deli$ in the above three expressions,
the argument resulting in (\ref{eng}) implies that in all three
cases, \be V^{(i)}_n(\alpha;\deli)\le
\frac{V^{(i)}(\alpha)}{M^{(i)}}+\epsilon_n^{(i)},\quad
\epsilon_n^{(i)}=\epsilon_n^{(i)}(\alpha)\to 0, \quad n\to \infty,
\label{VI}\eu uniformly for all $\alpha\in\BR.$ This completes the
proof. \hfill\qed

\section{The local limit theorem and completion of the proof}
Local limit theorems are viewed as the main ingredient of the
Khintchine method. Theorem~\ref{loc} below says that a local
limit theorem holds for all three types of  structures obeying the
conditions of our Theorem~\ref{main}.
 \begin{theorem}[Local limit theorem]\label{loc}.
Let $\deli$, $i=1,2,3$ denote the solutions to the equations
(\ref{e1})-((\ref{e3}) respectively and let  the random variables
$Z_n^{(i)}, \ n\ge 1$ be defined as in (\ref{zdef}), where the
random variables $Y_k$ have distributions given in the paragraph
following (\ref{F}). Assume that condition $(iii^\prime)$ of
Theorem~\ref{main} holds for $i=1,2,3$. Then,
$$
 \PP(Z_n^{(i)}=n)\sim\frac{1}{\sqrt{2\pi {\rm Var(Z_n^{(i)})}}}
\sim\frac{1}{\sqrt{2\pi
K_2^{(i)}}}\left(\deli\right)^{1+r/2},\quad
 n\rightarrow\infty,\quad i=1,2,3,
$$
with constants $K_2^{(i)}$ defined by
$$
K_2^{(1)}=A\Gamma(r+2)\zeta(r+1),
$$
$$
K_2^{(2)}=A(1-2^{-r})\Gamma(r+2)\zeta(r+1)
$$
and
$$
K_2^{(3)}=A\Gamma(r+2).
$$
\end{theorem}

\nin {\bf Proof}\  We will find asymptotics for $\BP(Z_n=n)$ as
$n\rightarrow\infty$ for the  three types of random structures.
Following the pattern of the Khintchine method (see e.g.
\cite{frgr1,GS}), we set $\delta=\delta_n$ in (\ref{Y}) and
(\ref{ch}) and define $\alpha_0=\alpha_0(n)$ to be
\ber\label{t0def}
\alpha_0=\delta_n^{\frac{r+2}{2(r+1)}}\log^2n.\ena Then we have
 \ber\label{decomp}
 \PP(Z_n=n)=\int_{-1/2}^{1/2}\phi_n(\alpha)
e^{-2\pi in\alpha}d\alpha=I_1+I_2,
 \la{probeq}
 \ena
 where
$I_1=I_1(n)$ and $I_2=I_2(n)$ are defined to be \be\label{I1def}
I_1=\int_{-\alpha_0}^{\alpha_0}\phi_n(\alpha)e^{-2\pi
in\alpha}d\alpha \eu and \be\label{I2def}
I_2=\int_{-1/2}^{-\alpha_0}\phi_n(\alpha)e^{-2\pi
in\alpha}d\alpha +\int_{\alpha_0}^{1/2}\phi_n(\alpha)e^{-2\pi
in\alpha}d\alpha. \eu

Defining $B_n$ and $T_n$ by \be
  B_n^2=
\frac{d^2}{d\delta^2}\left(\log f_n(e^{-\delta_n})\right),\quad
\quad T_n=-\frac{d^3}{d\delta^3}\left(\log
f_n(e^{-\delta_n})\right), \label{BT}
\end{equation}
for $n$ fixed we have the expansion
 \begin{eqnarray}\label{expand2}
 \non\phi_n(\alpha)e^{-2\pi in\alpha}
&=&
\exp{\left(2\pi i\alpha(EZ_n-n)-2\pi^2\alpha^2B_n^2+O(\alpha^3)T_n\right)}\nonumber\\
&=&\exp{\left(-2\pi^2\alpha^2B_n^2+O(\alpha^3) T_n\right)},\quad
\alpha\rightarrow0. \label{exp}
 \end{eqnarray}
It can be checked that $B_n^2={\rm Var} Z_n$ and
$T_n=\sum_{j=1}^n\E(Y_j-\E Y_j)^3,$ by the argument leading to
\refm[EZ].

Now  (\ref{prime}) in  Lemma~\ref{prodasymp} and \refm[BT] tell us
that for all structures considered \be\label{Basymp}
(B_n^2)^{(i)}\sim K_2^{(i)}(\delta_n^{(i)})^{-r-2}, \eu and
$$
T_n^{(i)}\sim K_3^{(i)}(\delta_n^{(i)})^{-r-3},
$$
where $K_2^{(i)},K_3^{(i)}>0$, $i=1,2,3$, are constants depending
on the type of the structure and while  $K_2^{(i)}, \ i=1,2,3$
are as in the statement of the theorem.

 Therefore, considering (\ref{t0def}), we find
that in all three cases, \be\label{big}
\lim_{n\rightarrow\infty}{B_n^2\alpha_0^2}=\infty {\rm \ \ and \ \
} \lim_{n\rightarrow\infty}{T_n\alpha_0^3}=0. \eu
  Combining (\ref{expand2}) with (\ref{big}), we deduce that
 \be\label{phiexp}
 \phi_n(\alpha)e^{-2\pi in\alpha}\sim\exp{\left(-2\pi^2\alpha^2B_n^2\right)},\quad
 n\rightarrow\infty,\quad |\alpha|\leq \alpha_0.
 \eu
  Finally, using (\ref{I1def}), (\ref{big}) and (\ref{phiexp})
gives us \be\label{I1asymp}
 I_1\sim\int_{-\alpha_0}^{\alpha_0}{\exp{\left(-2\pi^2\alpha^2
 B_n^2\right)}}d\alpha
 \sim (2\pi B_n)^{-1}
 \int_{-\infty}^{\infty}{e^{-\halfe
 \alpha^2}}d\alpha=\frac{1}{\sqrt{2\pi B_n^2}},\quad  n\rightarrow\infty.
\eu

The next step of the proof is to show that
$I_2=o(I_1)$, $n\rightarrow\infty$.
At this step condition $(iii^\prime)$ of Theorem~\ref{main} plays
a key role. Because of the asymptotic
$\sqrt{\deli}=o(\alpha_0^{(i)}), \ \ n\to \infty,\ i=1,2,3,$ we
can use condition $(iii')$ to bound the quantity $V^{(i)}(\alpha)$
defined in Lemma~\ref{gtau} by
$$
 V^{(i)}(\alpha)\le\left(1+{r\over 2}+\epsilon\right) M^{(i)}\log\deli,
 \quad \alpha_0\leq \vert\alpha\vert\leq 1/2, \ i=1,2,3.
$$
Hence,  Lemma~\ref{gtau} and the fact that in condition $(iii')$,
$\epsilon>0$   give

$$
|\phi_n^{(i)}(\alpha)|=o(\left(\deli\right)^{1+{r\over
2}})(1+\epsilon_n^{(i)}),\quad \alpha_0\leq |\alpha|\leq 1/2, \
n\to\infty, \quad i=1,2,3.
$$
From the definition (\ref{I2def}) and  the asymptotic
(\ref{Basymp}) we have \be\label{I2est}
I_2=o\left(\left(\deli\right)^{1+r/2}\right)=o(I_1). \eu Lastly
from (\ref{probeq}), (\ref{Basymp}), (\ref{I2est}) and
(\ref{I1asymp}), we derive the following asymptotic expression for
$\PP(Z_n^{(i)}=n), \ i=1,2,3$:
$$
 \PP(Z_n^{(i)}=n)\sim\frac{1}{\sqrt{2\pi (B_n^2)^{(i)}}}
\sim\frac{1}{\sqrt{2\pi K_2^{(i)}}}\left(\deli\right)^{1+r/2},\quad
 n\rightarrow\infty.
$$
\hfill\qed \vskip 1cm

To complete the proof of Theorem~\ref{main} it is left to
substitute the asymptotic expressions implied by our results for
the three factors  in the representation (\ref{rep}) when
$\delta=\deli$.

\section{Concluding remarks}

(i) Under the stated conditions on
parameters $b_k,$ the asymptotic formulae \refm[casympmul],
\refm[sl], \refm[as], for $c_n^{(i)}, \ i=1,2,3$ have a striking
similarity, all of them being of the form:
$$c_n\sim \chi_1 n^{\chi_2}\exp{\left(\chi_3n^{\frac{r}{r+1}}\right)}, \quad
n\to \infty,$$ where we have denoted by $\chi_1, \chi_2,\chi_3$
the constants that depend on the type of a structure and it
parameters. A simple analysis that takes into account that
$\zeta(r+1)>1$ reveals that, asymptotically in $n$,
$c_n^{(1)}>c_n^{(2)}$ and $c_n^{(1)}>c_n^{(3)},$ where the first
fact follows obviously from the definition of selections.

(ii)The following  observation is also in order. It turns out that
each one of the three combinatorial structures obeying the above
conditions behave very much alike the one with parameters
$b_k=k^{r-1}, \ r>0.$ According to the classification suggested in
\cite{bar} the latter structures are called expansive. In this
respect, combinatorial structures obeying Meinardus' conditions
(as well as their extensions as defined in the present paper) can
be viewed as quasi-expansive.

\ignore{
(iii) {\bf Discussion of condition $(iii)^\prime.$ of
Theorem~\ref{main}} Assume that
 a combinatorial structure obeys
Meinardus conditions (i) and (ii), while for all sufficiently
small $\delta>0$  and some $1/2<\alpha>\sqrt{\delta},$ \be
C_1\log(\delta^{-1})< 2\sum_{k=1}^\infty
 {b_k{e^{-k\delta}\sin^2(\pi k\alpha)}}< M^{(i)}(1+\frac{r}{2})\log (\delta^{-1}),
\label{iop}\end{equation} where $0<C_1< M^{(i)}(1+\frac{r}{2}),$
which says that in the case considered condition
$(iii)^\prime$ is violated. D: So what?? May be remove the remark?
}

(iii) We hope that the approach of this paper can be applied as
well to other enumeration problems, in particular to enumeration
of structures with constraints on the number of summands
(components) (see e.g. \cite{rom}).

 \vskip.5cm {\bf
Acknowledgement.} We are grateful to Prof. Harold Shapiro
 who pointed to us Wiener-Ikehara theorem, during his visit at
 Technion in 2006.

\end{document}